\documentclass[12pt,oneside]{amsart}

\usepackage{amscd}
\usepackage{amssymb}
\usepackage[latin1]{inputenc}
\usepackage[french]{babel}

\def\Z{\mathbb{Z}}

\def\ZZ{\Z \oplus \Z}
\def\R{\mathbb{R}}
\def\H{\mathbb{H}}

\def\P{\pi_{1}}

\def\G{\Gamma}

\newtheorem{thm}{Théorème}

\newtheorem{conj}{Conjecture}

\title{Sur la conjecture des fibrés de Seifert}

\begin{document}
\maketitle 
\begin{center}
{\sc Jean-Philippe PR\' EAUX}\footnote[1]{Centre de recherche de l'Armée de l'Air, Ecole de l'air, F-13661 Salon de
Provence air}\ \footnote[2]{Centre de Math\'ematiques et d'informatique, Universit\'e de Provence, 39 rue
F.Joliot-Curie, F-13453 marseille
cedex 13\\
\indent {\it E-mail :} \ preaux@cmi.univ-mrs.fr}
\end{center}
\begin{abstract}
Nous rappelons l'historique de la démonstration  de la conjecture des fibrés de Seifert ainsi que ses motivations et
ses diverses généralisations.
\end{abstract}

\section*{Introduction}
La conjecture des fibrés de Seifert a pour objet la caractérisation des fibrés de Seifert dans la classe des 3-variétés
orientables irréductibles à $\pi_1$ infini par une propriété de leur groupe fondamental : nommément l'existence d'un
sous-groupe normal cyclique non trivial. Il s'agit d'une question de grande importance en topologie des 3-variétés car
reliée à diverses autres conjectures de grand impact, comme nous ne le verrons par la suite. C'est désormais un
théorème dont la démonstration a necessité le travail de nombreux mathématiciens durant près d'un demi-siècle ; l'un de
ces théorèmes "monstrueux" qu'ont vu apparaître les mathématiques du vingtième siècle. Nous rappelons ici ses
motivations, ses diverses généralisations ainsi que l'historique de leur démonstration. Le lecteur est supposé être
familiarisé avec la topologie des 3-variétés.

\section{Enoncé de la conjecture des Fibrés de Seifert}

Un groupe fondamental infini d'un fibré de Seifert contient un sous-groupe normal cyclique infini. La conjecture des
fibrés de Seifert (ou CFS) caractérise les espaces fibrés de Seifert dans la classe des 3-variétés orientables
irréductibles à $\pi_1$ infini par le biais de cette propriété. Elle s'énonce :

\begin{conj}
Soit $M$ une 3-variété orientable irréductible dont le $\pi_1$ est
infini et contient un sous-groupe normal cyclique non trivial.
Alors $M$ est un fibré de Seifert.
\end{conj}

Elle se généralise au cas non-orientable :

\begin{conj}
Soit $M$ une 3-variété $\Bbb P_2$-irréductible dont le $\pi_1$ est
infini et contient un sous-groupe normal cyclique non trivial.
Alors $M$ est un fibré de Seifert.
\end{conj}

\noindent{\bf Remarques :} -- Des arguments classiques de topologie algébrique et le théorème de la sphère montrent
qu'une 3-variété $\Bbb P_2$-irréductible à $\pi_1$ infini a un groupe libre de torsion. Ainsi dans ces conjectures on
peut remplacer "... dont le $\pi_1$ est infini et contient un sous-groupe normal cyclique non trivial." par : "... dont
le $\pi_1$ contient un sous-groupe normal cyclique infini.".\\
-- Un fibré de Seifert orientable est soit irréductible soit homéomorphe à $S^1\times S^2$ ou à $\Bbb P^3\#\Bbb P^3$.
Comme conséquence du théorème de Kneser-Milnor, une 3-variété orientable $M$ non irréductible contenant un sous-groupe
normal cyclique non trivial est soit $S^1\times S^2$, soit $M'\# C$ avec $M'$ irréductible et $C$ simplement connexe,
soit a pour $\pi_1$  le groupe diédral infini. Si l'on accepte la conjecture de Poincaré (aujourd'hui démontrée par les
travaux de Perelman, récompensé d'une médaille Fields en 2006)  $M$ ne peut être qu'obtenu à partir de $S^1\times S^2$,
$\Bbb P^3\# \Bbb P^3$ ou d'une 3-variété irréductible en lui retirant un nombre fini
de boules. \\
-- Un fibré de Seifert non-orientable est quant à lui soit $\Bbb P^2$-irréductible soit $S^1\ltimes S^2$ ou $\Bbb
P^2\times S^1$. Une large classe de 3-variétés irréductibles et non $\Bbb P^2$-irréductible n'admettent aucune
fibration de Seifert tandis que leur $\pi_1$ contient un $\Z$ distingué. Nous verrons que le résultat se généralise
cependant aussi dans ce cas en considérant les fibrés de
Seifert mod $\Bbb P$.\medskip\\
\indent La CFS est aujourd'hui devenue un théorème, grâce au travail commun (énorme) des Mathématiciens : Waldhausen
(\cite{wald}), Gordon et Heil (\cite{gh}), Jaco et Shalen (\cite{js}), pour le cas Haken ; Scott (\cite{scott}), Mess
(\cite{mess}), Tukia (\cite{tukia}), Casson et Jungreis (\cite{cass}), Gabai (\cite{gabai}), pour le cas orientable non
Haken ; et Heil et Whitten (\cite{whit,{heil-whitten}}) pour le cas non-orientable ; on peut aussi citer  Maillot
(\cite{maill}, \cite{maillot}) et Bowditch (\cite{bow2}) qui donnent une preuve alternative incluant le résultat
non-publié de Mess.

Remarquons aussi que le théorème se généralise de diverses façons : aux 3-variétés ouvertes, aux 3-orbiétés
(\cite{maillot}), et aux groupes $PD(3)$ (\cite{bow2}), ou encore en affaiblissant la condition d'existence d'un $\Z$
distingué en l'existence d'une classe de conjugaison finie non triviale (\cite{hp}).

\section{Motivations}

 Historiquement, trois questions d'importance
ont motivé la CFS. D'abord la conjecture du centre (1960's), puis le strong torus theorem de Scott (1978), et enfin le
programme de géométrisation de Thurston (1980's).
\subsection{Conjecture du centre} Il s'agit du problème 3.5 de la
liste de Kirby, attribué à Thurston. :\smallskip\\
\noindent{\bf Conjecture : } {\sl Soit $M$ une 3-variété
orientable irréductible ayant un $\pi_1$ infini et un centre non
trivial, alors $M$ est un fibré de Seifert.}\smallskip

C'est clairement un corollaire immédiat de la conjecture des
fibrés de Seifert.
Après avoir été observé et successivement prouvé pour tous les compléments de noeuds (Murasugi (\cite{mu}), Neuwirth
(\cite{ne}) en 1961 pour les noeuds alternés, et Burde et Zieschang (\cite{bz}) en 1966 pour tous les noeuds), il a été
démontré en 1967 par Waldhausen (\cite{wald}) dans le cas où $M$ est Haken, qui est somme toute plus général que les
précédents.

Le sous-groupe normal cyclique infini d'un $\P$ de fibré de Seifert orientable est central si et seulement si la base
de la fibration associée est orientable. Un fibré de Seifert à $\pi_1$ infini et à base non-orientable n'a quant à lui
pas de centre, mais un sous-groupe normal cyclique non trivial. La CFS généralise la conjecture du centre en ce sens.

\subsection{Le théorème du tore} La conjecture appelée "théorème
du tore" s'énonce :\smallskip\\
{\bf Conjecture :} {\sl Soit $M$ une 3-variété orientable
irréductible avec $\pi_1(M)\supset \ZZ$. Alors soit $M$ contient
un tore incompressible, soit $M$ est un (small) fibré de
Seifert.}\smallskip

Elle a été démontrée dans le cas où $M$ est Haken par Waldhausen en 1968 (annoncé dans \cite{w2}, écrit et publié par
Feustel dans \cite{f1,f2}). Avec la théorie des 3-variétés {\sl suffisament grandes} dévellopée par Haken celà mena
(1979)
à l'apparition de la décomposition {\sl Jaco-Shalen-Johansen} des 3-variétés Haken.\smallskip\\
\indent
En 1978 Scott démontre (\cite{stt}) le "strong torus theorem" :\smallskip\\
{\bf Strong torus theorem : } {\sl Soit $M$ une 3-variété
orientable irréductible, avec $\pi_1(M)\supset \ZZ$. Alors soit
$M$ contient un tore incompressible, soit $\pi_1(M)$ contient un
sous-groupe normal cyclique non trivial.}\smallskip

Avec le strong torus theorem, pour prouver le théorème du tore il suffit de prouver la conjecture des fibrés de
Seifert.

\subsection{La conjecture de géométrisation}
La conjecture de géométrisation de Thur\-ston entraînerait la classification des 3-variétés.
Elle suppute que les pièces obtenues dans la décomposition topologique canonique d'une 3-variété orientable, le long de
sphères disques et tores essentiels, ont un intérieur métrisable au sens de Riemann de façon complète et localement
homogène. L'intérieur des pièces est alors modelé sur l'une des 8 géométries homogènes 3-dimensionnelles : les 3
isotropes, l'elliptique $\Bbb S^3$, l'euclidienne $\Bbb E^3$ et l'hyperbolique --la générique-- $\Bbb H^3$, les 2
géométries produits $\Bbb S^2\times\R$ et $\Bbb H^2\times \R$, et les 3 géométries twistées, $Nil$, $Sol$, et le
revêtement universel de $SL_2\R$.
Une 3-variété orientable admet une fibration de Seifert
exactement lorsque son intérieur est modelé sur l'une des 6 géométries autres que $Sol$ et $\H^3$.\\

Thurston ({\sl et al...}) a montré la conjecture de géométrisation dans le cas Haken :\smallskip\\
{\bf Théorème de géométrisation de Thurston : } {\sl Si $M$ est Haken
alors $M$ vérifie la conjecture de géométrisation de Thurston.}\\

Pour prouver la conjecture dans les cas restants il suffit de
prouver sous l'hypothèse que $M$ est une 3-variété orientable
irréductible fermée, les
3 conjectures :\smallskip\\
{\sl-1- Si $\pi_1(M)$ est fini alors $M$ est elliptique
(conjecture
d'orthogonalisation).\smallskip\\
-2- Si $\pi_1(M)$ est infini et contient un sous-groupe normal cyclique non trivial alors $M$ est un fibré de Seifert
(conjecture des fibrés
de Seifert.)\smallskip\\
-3- Si $\pi_1(M)$ est infini et ne contient pas de sous-groupe
normal cyclique non trivial, alors $M$ est géométrisable.} \smallskip\\
A la lumière du strong torus theorem, du théorème de géométrisation de Thurston et du fait qu'une 3-variété fermée
modelée sur l'une des 7 géométries non hyperboliques a un $\pi_1$ contenant $\ZZ$,
le point 3 devient :\smallskip\\
{\sl-3- Si $\pi_1(M)$ est infini et ne contient pas $\ZZ$ alors
$M$
est hyperbolique (conjecture d'hyperbolisation).}\smallskip\\

Ainsi la conjecture des fibrés de Seifert apparaît comme un des 3 pans de la conjecture de géométrisation... Celui
considéré comme le plus facile ; par ailleurs le premier à avoir été démontré. Les deux derniers pans sont quant à eux
 en passe d'être démontrés par les travaux de Perelman  ;
 il a suivi le programme d'Hamilton utilisant le flôt de Ricci.

\section{Historique de la démonstration de CFS}
\subsection{Le cas Haken orientable}L'historique de la preuve est le suivant
:\smallskip\\
{\bf 1967.} Waldhausen (\cite{wald}, \cite{f1,f2}) montre qu'une 3-variété Haken $M$ a un $\pi_1$ à centre non trivial
si et seulement si c'est un fibré de Seifert à base orientable.

Il motive ainsi la conjecture des fibrés de Seifert, et résout le
cas Haken lorsque le sous-groupe cyclique normal est
central.\smallskip\\
{\bf 1975.} Gordon et Heil (\cite{gh}) montrent partiellement la CFS dans le cas Haken : soit $M$ est un Seifert soit
obtenu en recollant deux copies d'un $I$-fibré non trivial sur une surface non-orientable. Ils réduisent ainsi  les cas
Haken restants à ces dernières 3-variétés.\smallskip\\
{\bf 1979.} Jaco et Shalen (\cite{js}), et indépendamment MacLachlan (non publié) achèvent la preuve pour les
3-variétés orientables Haken restantes.

\subsection{Le cas non Haken, orientable} Avec le cas Haken déjà établi, il est suffisant de se restreindre au
des 3-variétés closes.
 La preuve a procédé ainsi
:\medskip\\
 \noindent {\bf 1983.} Scott (\cite{scott})
(généralisant un résultat de
Waldhausen dans le cas Haken) montre que : \\
{\sl \indent Soient $M$ et $N$ deux 3-variétés closes orientables
irréductibles ; où $N$ est un fibré de Seifert à $\pi_1$ infini.
Si $\pi_1(M)$ et $\pi_1(N)$ sont isomorphes, alors $M$ et $N$ sont
homéomorphes.}\smallskip\\
\noindent {\bf Remarque :} Par hypothèse, avec le théorème de la
sphère et des arguments classiques de topologie algébrique, $M$ et
$N$ sont des $K(\pi,1)$. Avec le théorème de Moise (...), on peut
les considérer dans la catégorie PL. Ainsi la condition
"$\pi_1(M)$ et $\pi_1(N)$ sont isomorphes" peut-être remplacée par
"$M$ et $N$
ont même type d'homotopie".\smallskip\\
\indent
Ainsi Scott réduit la conjecture des fibrés de Seifert à :\smallskip\\
\indent {\sl Si $\G$ est le groupe d'une 3-variété close
irréductible orientable contenant un $\Z$ distingué, alors $\G$
est le groupe d'un fibré de Seifert clos orientable.}\smallskip

Il remarque en outre que $\G$ est le groupe d'un fibré de Seifert (clos orientable) si et seulement si $\G/\Z$ est le
groupe d'une 2-orbiété (close éventuellement non-orientable) ; c'est aussi le lemme 15.3 de \cite{bow}. Ainsi
finalement, il réduit la preuve de la conjecture des fibrés de Seifert à
prouver la conjecture :\smallskip\\
\indent{\sl\bf Si $\G$ est le groupe d'une 3-variété orientable
close contenant $\Z$ distingué, alors $\G/\Z$ est le groupe d'une
2-orbiété (close)}.\\

\noindent {\bf Fin 80'.} Mess (\cite{mess}) montre dans un papier
malheureusement jamais publié que :\smallskip\\
\indent{\sl Si $M$ est close orientable irréductible avec
$\G=\pi_1(M)$
contenant un sous-groupe cyclique infini $C$ :\\
-- Le revêtement de $M$ associé à $C$ est homéomorphe à la variété
ouverte $D^2\times S^1$.\\
-- L'action de rev\^etement de $\G/C$ sur $D^2\times S^1$ descend
en une action vaguement définie sur $D^2\times 1$ : dans sa
terminologie $\G/C$ est coarse quasi-isometric au plan euclidien
ou
hyperbolique.\\
-- Dans le cas où $\G/C$ est coarse quasi-isometric au plan
euclidien, alors c'est le groupe d'une 2-orbiété, et donc
$M$ est un fibré de Seifert (avec le résultat de Scott).\\
-- Dans le cas restant où $\G/C$ est coarse quasi-isometric au
plan hyperbolique, $\G/C$ induit une action sur le cercle à
l'infini, qui fait de $\G/C$ un groupe de convergence.}\smallskip\\
\indent
 Un {\sl groupe de convergence} est un groupe  $G$ agissant par homéomorphisme préser\-vant l'orientation sur le
cercle, de façon à ce que si l'on considère l'ensemble $T$ des triplés ordonnés : $(x,y,z)\in S^1\times S^1\times S^1$,
$x\not=y\not=z$, tels que $x,y,z$ apparaissent dans cet ordre sur $S^1$ en tournant dans le sens positif, l'action de
$G$ induite sur $T$ est libre et proprement
discontinue.\smallskip\\
Avec le travail de Mess, la preuve de la conjecture des fibrés de Seifert se
réduit à montrer la conjecture :\\
{\sl\bf \indent Les groupes de convergence sont des groupes de
2-orbiétés.}\\

\noindent {\bf 1988.} Tukia montre (\cite{tukia}) que certains (lorsque $T/\G$ non compact et $\G$ n'a pas d'éléments
de torsion d'ordre $>3$) groupes de convergence sont des groupes Fuchsiens, et en particulier des
$\pi_1$ de 2-orbiétés.\\

\noindent {\bf 1992.} Gabai montre (\cite{gabai}) indépendamment des autres travaux et en toute généralité que les
groupes de convergence sont des groupes fuchsiens (agissant sur $S^1$, à conjugaison dans $Homeo(S^1)$ près, comme la
restriction à $S^1=\partial \Bbb H^2$
de l'action d'un groupe fuchsien sur $\bar{\Bbb H}^2$).\\

\noindent {\bf 1994.} Parallèlement Casson et Jungreis montrent les cas laissés
restants par Tukia (\cite{cass}).\\

{\bf Ainsi la conjecture des fibrés de Seifert est prouvée, et
avec elle la conjecture du centre, le théorème du tore, et une des
trois conjectures de géométrisation.}\\

\noindent{\bf 1999.} Bowditch obtient une preuve différente de la conjecture des fibrés de Seifert (\cite{bow},
\cite{bow2}). Sa preuve se généralise à d'autres types
de groupes, comme les groupes $PD(3)$ par exemple.\\

\noindent{\bf 2000.} Maillot, dans sa thèse, étend les techniques de Mess pour établir une preuve de la conjecture des
fibrés de Seifert dans le cas des 3-variétés ouvertes ainsi que des 3-orbiétés, généralisant le cas des 3-variétés. Ce
résultat, déjà connu comme conséquence de la conjecture des fibrés de Seifert (en tant que prouvée par Mess, {\sl et
al}...) et du théorème des orbiétés de Thurston, a le mérite d'être redémontré en n'utilisant aucun de ces résultats,
et d'obtenir comme corollaire une preuve propre de la conjecture de Seifert et des techniques de Mess.

L'argument de Mess est repris dans \cite{maillot} où il se ramène
aux groupes quasi-isométriques à un plan riemannien, et il montre
dans \cite{maill} que ces derniers sont virtuellement des groupes
de surface.

\subsection{Le cas non-orientable}
Le cas non-orientable a été traité par Whitten et Heil.\\

\noindent{\bf 1992.} Whitten montre dans \cite{whit} la CFS pour $M$ non-orientable irréductible, qui n'est pas un faux
$\Bbb P_2\times S^1$, et tel que $\pi_1(M)$ ne contient pas $\Z_2*\Z_2$. Il obtient en particulier la conjecture 2, ou
CFS dans le cas $\Bbb
P_2$-irréductible.\\

\noindent{\bf 1994.} Heil et Whitten dans \cite{heil-whitten} caractérisent les 3-variétés irréductibles
non-orientables ne contenant pas un faux $\Bbb P_2\times I$ et dont le $\pi_1$ contient  un $\Z$ distingué et
éventuellement $\Z_2*\Z_2$ : ce sont les fibrés de Seifert mod $\Bbb P$ : obtenus à partir d'une 3-orbiété fibré de
Seifert  en retirant tous les voisinages des points singuliers homéomorphes à des cônes sur $\Bbb P_2$ : leur bord
contient en général des $\Bbb P_2$ et dans le cas contraire il s'agit d'un fibré de Seifert. Le résultat devient :

\begin{thm} Si $M$ est une 3-variété non-orientable irréductible ne contenant pas un faux $\Bbb P_2\times I$ avec $\pi_1(M)$
contenant un sous-groupe cyclique non trivial alors $M$ est soit $\Bbb P_2\times I$ soit un fibré de Seifert mod $\Bbb
P$.
\end{thm}

Ils montrent aussi le théorème du tore dans le cas non-orientable, et déduisent la géométrisation dans le cas non
orientable ({\sl i.e.} le revêtement d'orientation est géométrisable) si la 3-variété ne contient pas de faux $\Bbb
P_2\times I$.

\subsection{Avec la conjecture de Poincaré}\hfill\\

\noindent {\bf 2003-06.} Les travaux de Perelman (2003) récompensés en 2006 par une médaille Fiels établissent la
conjecture de Poincaré. Ainsi toute 3-variété fermée simplement connexe est une 3-sphère et il n'existe ni fausse boule
ni faux $\Bbb P_2\times I$. Aussi on peut supprimer dans le théorème de CFS l'hypothèse d'irréducibilité sans trop
alourdir l'énoncé.
\begin{thm} Soit $M$ une 3-variété dont le $\pi_1$ est infini et contient un sous-groupe normal cyclique non
trivial. Alors après avoir bouché les sphères de $\partial M$ avec des boules on obtient soit une somme connexe de
$\Bbb P_2\times I$ avec lui-même ou avec $\Bbb P_3$, soit un fibré de Seifert mod $\Bbb P$.
\end{thm}
On obtient un fibré de Seifert exactement lorsque $\P(M)$ ne contient pas $\Z_2$ ou de façon équivalente lorsque
$\partial M$ ne contient pas de $\Bbb P_2$.

\subsection{Groupes de 3-variétés à classes de conjugaison infinies}\hfill\\

\noindent {\bf 2005.} de la Harpe et l'auteur montrent dans \cite{hp} que dans le cas des 3-variétés et des groupes
$PD(3)$ l'hypothèse "{\sl contient un sous-groupe normal cyclique non trivial}" peut être affaiblie en "{\sl contient
une classe de conjugaison finie non triviale}" ou de façon équivalente en "{\sl a une algèbre de Von-Neumann qui n'est
pas un facteur de type $II-1$}".

\end{document}